\documentclass[12pt]{amsart}
\usepackage{amssymb,amsmath}

\headsep=1cm \numberwithin{equation}{section}
\newtheorem{theorem}{Theorem}[section]
\newtheorem{lemma}[theorem]{Lemma}

\newtheorem{corollary}[theorem]{Corollary}

\theoremstyle{definition}

\theoremstyle{remark}

\newcommand{\im}{\operatorname{im}}

\newcommand{\Ext}{\operatorname{Ext}}

\newcommand{\fm}{\frak{m}}

\newcommand{\fa}{\frak{a}}

\begin{document}
\author[Mafi ]{Amir Mafi}
\title[Some properties of generalized local cohomology modules]
{Some properties of generalized local cohomology modules}

\address{A. Mafi,  Institute of Mathematics, University for Teacher
Education, 599 Taleghani Avenue, Tehran 15614, Iran and Department
of Mathematic, Arak University, Arak, Iran.}
\email{a-mafi@araku.ac.ir}



\begin{abstract} Let $R$ be a commutative Noetherian ring, $\fa$ an ideal of
$R$, $M$ and $N$ be two finitely generated $R$-modules. Let $t$ be
a positive integer. We prove that if $R$ is local with maximal
ideal $\fm$ and
 $ M\otimes_R N$ is of finite length then $H_{\fm}^t(M,N)$ is of
 finite length for all $t\geq 0$ and $l_R(H_{\fm}^t (M,N))\leq
 \sum _{i=0}^t l_R(\Ext_R^i(M,H_{\fm}^{t-i}(N)))$.
This yields, $l_R(H_{\fm}^t(M,N))=l_R(\Ext_R^t(M,N))$.\\
Additionally, we show that
  $\Ext_R^i(R/{\fa},N)$ is Artinian for all $ i\leq t$ if and only
  if $H_{\fa}^i(M,N)$ is Artinian for all $i\leq t$. Moreover, we show
that whenever $\dim (R/{\fa})=0$ then $H_{\fa}^t(M,N)$ is Artinian
for all $t \geq 0$.
\end{abstract}

\maketitle

\section{Introduction}

Throughout this paper, $R$ is a commutative Noetherian ring with
identity and $\fa$ is an ideal of $R$. For an $R$-module $N$ and
non-negative integer $t$, the  local cohomology module
$H_{\fa}^t(N)$ was first introduced and studied by Grothendieck
[{\bf 4}]. For example he showed that $H_{\fm}^t(N)$ is Artinian
for all $t$, whenever $N$ is finitely generated and $R$ is local
 with maximal ideal ${\fm}$. One of the general problem in local
cohomology is to find when the local cohomology module
$H_{\fa}^t(N)$ is Artinian (see [{\bf 6}, Problem 3]). In [{\bf
7}], Melkersson proved that for an $R$-module $N$, $H_{\fa}^i(N)$
is Artinian for all $i\leq t$ if and only if $\Ext_R^i(R/{\fa},N)$
is Artinian for all $i\leq t$. The latter result leads us to
consider the same result for generalized local cohomology modules
which was introduced and studied by Herzog [{\bf 5}] (see also
[{\bf 10}]). For each $i\geq 0$, the generalized local cohomology
functor $H_{\fa}^i(., .)$ is defined by
$$H_{\fa}^i(M,N)=\underset{n}{\varinjlim}\Ext_R^i(M/\fa^nM,N)$$ for
all $R$-modules $M$ and $N$. Clearly, this is a generalization of
the usual local cohomology module. The reader is referred to
articles [{\bf 2}], [{\bf 1}]and [{\bf 3}] for more results on
generalized local cohomology. Our main purpose in this paper is to
show the following.

\begin{theorem} Let ${\fa}$ be an ideal of $R$, $M$ a finitely generated $R$-module,
and let $t$ be a positive integer.\\
($\alpha$) If $R$ is local with maximal ideal $\fm$ and $N$ is a
finitely generated $R$-module such that $M\otimes_R N$ is of finite
length, then  $l_R(H_{\fm}^t(M,N))\leq \sum_{i=0}^t
l_R(\Ext_R^i(M,H_{\fm}^{t-i}(N)))$ for any $t\geq 0$.\\
 ($\beta$) Let $N$ be an $R$- module. Then the following statements are
equivalent:\\
(i) $\Ext_R^i(R/{\fa},N)$ is an Artinian $R$-module for all  $i\leq t$.\\
(ii) $H_{\fa}^i(M,N)$ is an Artinian $R$-module for all $i\leq t$.
\end{theorem}

Clearly ($\alpha$) extends the result [{\bf 9}, Theorem 3.2] and
($\beta$) is equivalent with the main result [{\bf 7}, Theorem
1.2].

\section{The results}

We start this section with the following lemma.

\begin{lemma} Let $R$ be a local ring with maximal ideal $\fm$, $M$
and $N$ be two finitely generated $R$-modules such that
$M\otimes_RN$ is of finite length. Then $H_{\fm}^t(M,N)$ is of
finite length for all $t\geq 0$.
\end{lemma}

{\bf Proof.} Consider the Grothendieck spectral sequence [{\bf 8},
Theorem 11.38]
$$E_2^{p,q}:=\Ext_R^p(M,H^q_{\fm}(N))\underset{p}{\Longrightarrow}
H_{\fm}^{p+q}(M,N).$$ for $t\geq 0$, we have a finite filtration
$$0=\phi^{t+1}H^t\subseteq\phi^tH^t\subseteq\dots\subseteq\phi^1H^t\subseteq\phi^
0H^t=H_{\fm}^t(M,N),$$ such that
$E_{\infty}^{i,t-i}=\phi^iH^t/\phi^{i+1}H^t$ for all $0\leq i\leq
t$. Now, since $ E_{\infty}^{i,j}$ is a homomorphic image of
$E_2^{i,j}$ for all $i,j\geq 0$, by [{\bf 9}, Lemma 3.1],
$E_2^{i,j}$ is of finite length for all $i,j\geq 0$. It therefore
follows that $E_\infty^{i,j}$ is of finite length for all $i,j\geq
0$ and that $\phi^{t-i}H^t$ is of finite length for all $0\leq
i\leq t$. Hence, using the exact sequence $$ 0\longrightarrow
\phi^1H^t\longrightarrow H_{\fm}^t(M,N)\longrightarrow
E_\infty^{0,t}\longrightarrow 0 (0\leq i\leq t)$$ we get
$H_{\fm}^t(M,N)$ is of finite length for all $t\geq 0$. $\Box$\\

The following theorem extends [{\bf 9}, Theorem 3.2].\\

 \begin{theorem} Let $R$ be a local ring with maximal
ideal $m$ and let $M$ and $N$ be two finitely generated
$R$-modules such that $M\otimes_RN$ is of finite length. Then
$l_R(H_{\fm}^t(M,N))\leq\sum_{i=0}^t
l_R(\Ext_R^i(M,H_{\fm}^{t-i}(N)))$ for all $t\geq 0$.
Consequently, $l_R(H_{\fm}^t(M,N))=l_R(\Ext_R^t(M,N))$.
\end{theorem}

{\bf Proof.} Let $2\leq i\leq t$. With the notation of [{\bf 8},
\S 11] we consider the exact sequences
$$ 0\longrightarrow \ker d_i^{i,t-i}\longrightarrow
E_i^{i,t-i}\longrightarrow E_i^{2i,t-2i+1}$$ and note that
$E_i^{i,t-i}=\ker d_{i-1}^{i,t-i}/\im d_{i-1}^{1,t-2}$ and
$E_i^{i,j}=0$ for all $j<0$. So, we have \\ $\ker
d_{t+2}^{i,t-i}\cong E_{t+2}^{i,t-i}\cong\ldots \cong
E_\infty^{i,t-i}$ for all $(0\leq i\leq t)$.
 Now, using the exact sequences
$$0\longrightarrow \phi^{i+1}H^t\longrightarrow \phi^iH^t
\longrightarrow E_\infty^{i,t-i}\longrightarrow 0 (0\leq i\leq
t)$$ and an argument similar to that used in 1.1 together with the
facts $E_\infty^{i,t-i}\cong\ker d_{t+2}^{i,t-i}\subseteq \ker
d_2^{i,t-i}\subseteq E_2^{i,t-i}$ for all $0\leq i\leq t$, one can
deduce that $l_R(H_{\fm}^t(M,N))\leq\sum_{i=0}^t
l_R(\Ext_R^i(M,H_{\fm}^{t-i}(N)))$. Moreover, there is a spectral
sequence $$
E_{2}^{i,j}=\Ext_R^i(M,H_{\fm}^{j}(N))\underset{i}{\Longrightarrow}
\Ext_{R}^{i+j}(M,N).$$ Hence, we have a finite filtration
$$0=\psi^{t+1}H^t\subseteq\psi^tH^t\subseteq\dots\subseteq\psi^1H^t\subseteq\psi^
0H^t=\Ext_{R}^{t}(M,N),$$ such that
$E_{\infty}^{i,t-i}=\psi^{i}H^{t}/{\psi^{i+1}H^{t}}$ for all
$0\leq i\leq t$. Now, using the same arguments as above, we get
$\psi^{i}H^{t}/{\psi^{i+1}H^{t}}=\varphi^{i}H^{t}/{\varphi^{i+1}H^{t}}$
for all $0\leq i\leq t$ and so the result follows. $\Box$\\

The following corollary extends [{\bf 9}, Corollary 3.3].

\begin{corollary} Let the situation be as in Theorem 2.2. Assume
that $N$ is Cohen-Macaulay with $\dim N=d$. Then
$H_{\fm}^t(M,N)\cong \Ext_R^t(M,N)$.
\end{corollary}

{\bf Proof.} By the same arguments as in the proof of Theorem 2.2
and [{\bf 9}, Corollary 3.3], we have $\Ext_R^t(M,N)\cong
\Ext_R^{t-d}(M,H_{\fm}^d(N))\cong E_\infty^{t-d,d}$, $$
0=\phi^tH^t=\phi^{t-1}H^t=\ldots =\phi^{t-d+1}H^t$$ and, also,
$$ \phi^ {t-d}H^t=\ldots=\phi^0H^t=H_{\fm}^t(M,N)$$ for all
$t$. It therefore follows that $H_{\fm}^t(M,N)\cong
\Ext_R^t(M,N)$. $\Box$\\

 The following theorem is related to [{\bf 7}, Theorem 1.2].

\begin{theorem} Let $N$ be an $R$-module and $t$ a positive integer.
Then the following conditions are equivalent:\\
(i) $\Ext_R^i(R/{\fa},N)$ is Artinian for all $i\leq t$.\\
(ii) $H_{\fa}^i(M,N)$ is Artinian for any finitely generated
$R$-module $M$ and all $i\leq t$.
\end{theorem}

{\bf Proof.} $(ii)\Longrightarrow (i)$ is immediate  by [{\bf
7}, Theorem 1.2 ].\\
 $(i)\Longrightarrow(ii).$ By [{\bf 7}, Theorem 1.2] $H_{\fa}^i(N)$ is Artinian for
all $i\leq t .$ Now by similar arguments as in the proof of 1.1.
One can see that $ E_\infty^{i,t-i}$ and $\phi^iH^t$ are Artinian
for all
 $0\leq i \leq t.$ We consider the exact sequences $$
0\longrightarrow \phi^ 1H^i\longrightarrow H_{\fa}^i(M,N)
\longrightarrow E_\infty^{0,i}\longrightarrow 0  (0\leq i\leq
t).$$ Hence $H_{\fa}^i(M,N)$ is Artinian for all $i\leq t$. $\Box$

\begin{corollary} Let the situation be as in Theorem 2.4. The
following conditions are equivalent:\\
(i) $\Ext_R^t(R/{\fa},N)$ is Artinian for all $t$.\\
(ii) $H_{\fa}^t(M,N)$ is Artinian for all finitely generated
$R$-modules $M$ and all $t$.\\
\end{corollary}

 The following corollary is a generalization of [{\bf 3}, Theorem 2.2].\\

 \begin{corollary} Let the situation be as in Corollary 2.5 and assume
 that $\dim(R/{\fa})=0$. Then $H_{\fa}^t(M,N)$ is Artinian for all
 finitely generated $R$-modules $M$ and all $t$.
 \end{corollary}

%%%%%%%%%%%%%%%%%%%%%%%%%%%%%%%%%%%%%%%%%%%%%%%%%%%%%%%%%%%%%%%%%%%%%%%%%%%%%

\end{document}